\documentclass{amsart}
\usepackage{amsmath}
\usepackage{latexsym,ulem}
\usepackage{amssymb}
\usepackage{graphics}

\input epsf%
\newtheorem{theorem}{Theorem}[section]

\newtheorem{lemma}{Lemma}[section]

\newtheorem{Definition}{Definition}[section]
\def\qed{{\hfill{\vrule height4pt width3pt depth2pt}}}

   \long\def\comment#1{}

\def\ad#1{\begin{aligned}#1\end{aligned}}  \def\b#1{\mathbf{#1}}
\def\a#1{\begin{align*}#1\end{align*}} \def\an#1{\begin{align}#1\end{align}}
\def\e#1{\begin{equation}#1\end{equation}} \def\d{\operatorname{div}}
\def\p#1{\begin{pmatrix}#1\end{pmatrix}} 
  \numberwithin{equation}{section}
\numberwithin{table}{section}
\numberwithin{figure}{section}

\def\boxit#1{\vbox{\hrule height1pt \hbox{\vrule width1pt\kern1pt
     #1\kern1pt\vrule width1pt}\hrule height1pt }}
 \def\lab#1{\boxit{\small #1}\label{#1}}
  \def\mref#1{\boxit{\small #1}\ref{#1}}
 \def\meqref#1{\boxit{\small #1}\eqref{#1}}

  \def\lab#1{\label{#1}} \def\mref#1{\ref{#1}} \def\meqref#1{\eqref{#1}}
\begin{document}

\newcommand{\disp}{\displaystyle}
\newcommand{\eps}{\varepsilon}
\newcommand{\To}{\longrightarrow}
\newcommand{\C}{\mathcal{C}}
\newcommand{\K}{\mathcal{K}}
\newcommand{\T}{\mathcal{T}}
\newcommand{\bq}{\begin{equation}}
\newcommand{\eq}{\end{equation}}
\long\def\comments#1{ #1}
\comments{ }
%%%%%%%%%%%%%%%%%%%%%%%%%%%%%%%%%%%%%%%%%%%%%%%%%%%%%%%%%%%%%%%%%%%%%%%%%%

\title[Conforming symmetric finite elements]{
A family of symmetric mixed finite
   elements for linear elasticity on tetrahedral grids
    }
\date{}

\author {Jun Hu}
\address{LMAM and School of Mathematical Sciences, Peking University,
  Beijing 100871, P. R. China.  hujun@math.pku.edu.cn}

 \author {Shangyou Zhang}
\address{Department of Mathematical Sciences, University of Delaware,
    Newark, DE 19716, USA.  szhang@udel.edu }

\begin{abstract}
A family of stable mixed finite elements for the linear   elasticity on
   tetrahedral grids are constructed,  where the
   stress is approximated by symmetric $H(\d)$-$P_k$ polynomial tensors
   and the displacement is approximated by
     $C^{-1}$-$P_{k-1}$ polynomial vectors,
    for all $k\ge 4$.
Numerical tests are provided.

  \vskip 15pt


 \vskip 15pt


\end{abstract}
\maketitle

\section{Introduction}

In the Hellinger-Reissner
mixed formulation of the linear elasticity equations,
   the stress is sought in $H(\d,\Omega,\mathbb {S})$
   and the displacement in $L^2(\Omega, \mathbb{R}^3)$.
It is a challenge to design stable mixed finites mainly due to
    the symmetric constraint of the stress  tensor $\mathbb {S}$.
To overcome this  difficulty,  ealier works adopted composite element
   techniques or weakly symmetric methods,  cf.  \cite{Amara-Thomas,
   Arnold-Brezzi-Douglas, Arnold-Douglas-Gupta,  Johnson-Mercier,
   Morley, Stenberg-1, Stenberg-2, Stenberg-3}.
In \cite{Arnold-Winther-conforming},  Arnold and Winther
   designed the first family of
    mixed finite element methods in 2D,
     based on  polynomial shape  function spaces.
   From then on,
      various  stable  mixed elements have been constructed,
        see \cite{Adams-Cockburn,Arnold-Awanou,Arnold-Awanou-Winther,
    Arnold-Winther-conforming,Awanou, Chen-Wang,
     Arnold-Winther-n, Gopalakrishnan-Guzman-n,Hu-Shi,
    Man-Hu-Shi, Yi-3D, Yi, Arnold-Falk-Winther,
    Boffi-Brezzi-Fortin, Cockburn-Gopalakrishnan-Guzman,
    Gopalakrishnan-Guzman, Guzman,Hu-Man-Zhang2014,Hu-Man-Zhang2013}.

As the displacement function is in $L^2(\Omega, \mathbb{R}^3)$,
  a natural discretization
   is the piecewise $P_{k-1}$ polynomial without interelement continuity.
It is  a long-standing and challenging problem
   if the  stress tensor can be discretized by
    an appropriate $P_k$ finite element subspace of
  $H(\d,\Omega,\mathbb {S})$.
Adams and Cockburn constructed such a mixed finite element in \cite{Adams-Cockburn}
    where the discrete stress space is the space
    of $H(\d,\Omega,\mathbb {S})$-$P_{k+2}$ tenors
    whose divergence is a $P_{k-1}$ polynomial on each tetrahedron,  for $k=2$.
The method was modified and extended to a family of elements, $k\ge 2$,
   by Arnold, Awanou and Winther \cite{Arnold-Awanou-Winther}.
Mathematically speaking,  these methods are two-order suboptimal.
In this paper,  we solve this open problem by constructing a
    suitable $H(\d,\Omega,\mathbb {S})$-$P_k$, instead of above $P_{k+2}$,
    finite element space for the stress discretization, for $k\ge 4$.
In these elements,  the symmetric stress tensor is approximated by
   the full $C^0$-$P_k$ space enriched by some so-called $H(\d)$
   edge-bubble functions locally on each tetrahedron.
A new way of proof  is developed to establish the stability of the
   mixed elements, by
     characterizing the divergence of local stress space.
This space of divergence of local stress space is
    exactly the subspace of $P_k$ displacements orthogonal to the local
    rigid-motion.
The optimal order error estimate is proved,
  verified by numerical tests of $P_4$ and $P_5$ mixed elements.

The rest of the paper is organized as follows. In Section 2, we
define the weak problem and the finite element method.
 In section 3,  we prove the well-posedness of the finite
    element problem, i.e. the discrete coerciveness and the
    discrete inf-sup condition.
  By which,  the optimal order convergence of the
   new element follows.
In Section 4, we provide some numerical results,
     using $P_4$ and $P_5$ finite elements.

\section{The family of finite elements}

Based on the Hellinger-Reissner principle, the  linear elasticity
     problem within a stress-displacement ($\sigma$-$u$) form reads:
Find $(\sigma,u)\in\Sigma\times V :=H({\rm div},\Omega,
    \mathbb {S}=\hbox{symmetric } \mathbb{R}^{3\times 3})
        \times L^2(\Omega,\mathbb{R}^3)$, such that
\an{\left\{ \ad{
  (A\sigma,\tau)+({\rm div}\tau,u)&= 0 && \hbox{for all \ } \tau\in\Sigma,\\
   ({\rm div}\sigma,v)&= (f,v) &\qquad& \hbox{for all \ } v\in V. }
   \right.\lab{eqn1}
}
Here the symmetric tensor space for stress $\Sigma$  and the
   space for vector displacement  $V$ are, respectively,
  \an{   \lab{S}
  H({\rm div},\Omega,\mathbb {S})
    &:= \Big\{ \sigma=\p{\sigma_{11} & \sigma_{12} &\sigma_{13} \\
        \sigma_{21}  & \sigma_{22} & \sigma_{23} \\
        \sigma_{31} &\sigma_{32}   & \sigma_{33}  }
     \in H(\d, \Omega)
    \ \Big| \ \sigma ^T = \sigma,  \Big\}, \\
     \lab{V}
     L^2(\Omega,\mathbb{R}^3) &:=
     \Big\{ \p{u_1 &  u_2 &u_3}^T
          \ \Big| \ u_i \in L^2(\Omega) \Big\}  .}
	This paper denotes by $H^k(T,X)$ the Sobolev space consisting of
functions with domain $T\subset\mathbb{R}^3$, taking values in the
finite-dimensional vector space $X$, and with all derivatives of
order at most $k$ square-integrable. For our purposes, the range
space $X$ will be either $\mathbb{S},$ $\mathbb{R}^3,$ or
$\mathbb{R}$.
$\|\cdot\|_{k,T}$ is the norm of $H^k(T)$. $\mathbb{S}$ denotes
the space of symmetric tensors, $H({\rm div},T,\mathbb{S})$
consists of square-integrable symmetric matrix fields with
square-integrable divergence. The H(div) norm is defined by
$$\|\tau\|_{H({\rm div},T)}^2:=\|\tau\|_{L^2(T)}^2
   +\|{\rm div}\tau\|_{L^2(T)}^2.$$
$L^2(T,\mathbb{R}^3)$ is the space of vector-valued functions
which are square-integrable.
Here, the compliance tensor
$A=A(x):\mathbb{S}~\rightarrow~\mathbb{S}$, characterizing the
properties of the material, is bounded and symmetric positive
definite uniformly for $x\in\Omega$.

This paper deals with a pure displacement problem \meqref{eqn1} with the
  homogeneous boundary condition that $u\equiv 0$ on
  $\partial\Omega$.
But the method and the analysis work for mixed boundary value problems
   and the pure traction problem.

The  domain $\Omega$ is subdivided by a family of quasi-uniform
  tetrahedral grids  $\mathcal{T}_h$ (with the grid size $h$).
We introduce the finite element space of order $k$ ($k\ge4$)  on $
   \mathcal{T}_h$.
The displacement space is the full $C^{-1}$-$P_{k-1}$ space
 \an{ \lab{Vh}
   V_h = \{v\in L^2(\Omega,\mathbb{R}^3) \ | \
         v|_K\in P_{k-1}(K, \mathbb{R}^3)\ \hbox{ for all } K\in\mathcal{T}_h \}.
    }

The discrete stress space of order $k$ ($k\ge 4$) is defined abstractly as
\an{\lab{Sh}
 \Sigma_h =\Big\{ \sigma
    &\in H(\d,\Omega,\mathbb{S}) \ \Big| \ \sigma|_K\in P_k(K, \mathbb{S})
     \   \forall K\in\mathcal{T}_h, \\
   \nonumber  & \ \sigma |_K(v_i)=\sigma |_{K'}(v_i) \quad
      \forall v_i\in \mathcal{V}_h \ \hbox{and } v_i\in K \cap K' \Big\}, }
where $\mathcal{V}_h$ is the set of vertices of the tetrahedral
          grid $\mathcal{T}_h$,  and $v_i$ is a common vertex of
          tetrahedra $K$ and $K'$.
 Computationally, for building a basis for $\Sigma_h$,
    we need to give another definition of $\Sigma_h$.
$\Sigma_h$ is a $H(\d )$ bubble enrichment of the $H^1$ space
  \an{\lab{tSh}
 \widetilde \Sigma_h =\Big\{~\sigma
    \in H^1(\Omega,\mathbb{S})  \ \Big| \ \sigma|_K\in P_k(K, \mathbb{S})
	        \ \forall K\in\mathcal{T}_h \Big\}.
            }

In computation, we still uses $6\times \dim P_k$ Lagrange nodal
   basis locally on each tetrahedron $K$, i.e., the standard basis
    for $H^1$ finite element space $\widetilde \Sigma_h$.
But globally, roughly speaking, we
   break each of $(k-1)=\dim P_{{k-2},1D}$ zero-flux (on all six edges)
      edge-bubble functions into $n_0$
    basis functions, where $n_0$ tetrahedra share this common edge,
   cf. Figure \mref{e-b},
   and break each of $(k-2)(k-1)/2=\dim P_{{k-3},2D}$
    zero-flux (on all four face triangles)
      facee-bubble functions into $2$
    basis functions, on the two tetrahedra sharing a common face triangle.
Here, on each triangle,  we have three sets of non-zero edge-bubble functions
   enriched,  all of which have a zero-flux on the triangle.
To avoid too much technical details,  we only
   define the local edge-bubble functions,  but we do not
   discuss on eliminating linearly dependent bubbles
    (with $H^1$-$P_k$ basis functions).

 \begin{figure}[htb]\setlength\unitlength{2pt}

 \begin{center}\begin{picture}(  100.,  60)(  0.,  -22)
     \def\lb{\circle*{0.8}}\def\lc{\vrule width1.2pt height1.2pt}

     \def\la{\circle*{0.3}}
 \put( -6.,  0.){$\b{x}_1$                                         }
 \put(  56.,  36.5999985){$\b{x}_2$                                         }
 \put(  101.,  0.){$\b{x}_3$                                         }
 \put(  43.0999985,  14.6999998){
 $\b{x}_0$                                         }
 \multiput(  33.30,  16.60)(  -0.096,   0.115){ 52}{\la} %arrow
 \multiput(  29.62,  19.80)(  -0.064,   0.136){ 20}{\la} %arrow
 \multiput(  30.82,  20.79)(  -0.122,   0.087){ 20}{\la} %arrow
 \multiput(  46.60,   5.00)(  -0.033,  -0.146){ 61}{\la} %arrow
 \multiput(  46.01,  -1.24)(  -0.068,  -0.134){ 20}{\la} %arrow
 \multiput(  44.49,  -0.90)(   0.005,  -0.150){ 20}{\la} %arrow
 \put(  25.,  23.6000004){$\b{n}_1$  }
 \put(  40., -7.){$\b{n}_2$                                         }
 \multiput(   0.00,   0.00)(   0.250,   0.000){400}{\la}
 \multiput(   0.00,   0.00)(   0.234,   0.088){170}{\la}
 \multiput(   0.00,   0.00)(   0.216,   0.126){277}{\la}
 \multiput( 100.00,   0.00)(  -0.243,   0.061){247}{\la}
 \multiput(  40.00,  15.00)(   0.177,   0.177){113}{\la}
 \multiput( 100.00,   0.00)(  -0.188,   0.165){212}{\la}
 \multiput(  20.00,  -6.00)(   0.000,   0.150){ 89}{\la} %arrow
 \multiput(  19.22,   4.50)(   0.038,   0.145){ 20}{\la} %arrow
 \multiput(  20.78,   4.50)(  -0.038,   0.145){ 20}{\la} %arrow

 \multiput(   5.00, -16.00)(   0.000,   0.150){119}{\la} %arrow
 \multiput(   4.22,  -1.13)(   0.038,   0.145){ 20}{\la} %arrow
 \multiput(   5.78,  -1.13)(  -0.038,   0.145){ 20}{\la} %arrow
  \multiput(  90.00,  30.00)(  -0.130,  -0.075){179}{\la} %arrow
 \multiput(  69.59,  17.41)(  -0.145,  -0.039){ 20}{\la} %arrow
 \multiput(  68.82,  18.77)(  -0.107,  -0.105){ 20}{\la} %arrow
 \put(  92.,  30.){$\lambda_1=0$                                     }
 \put(  91., -9.){$\lambda_0=0$                                     }
 \multiput(  65.00,   0.00)(   0.237,  -0.079){101}{\la}
 \multiput(  57.50,   2.50)(  -0.142,   0.047){ 21}{\la} %arrow
 \multiput(  57.10,   1.81)(  -0.126,   0.082){ 20}{\la} %arrow
 \multiput(  57.59,   3.29)(  -0.150,   0.010){ 20}{\la} %arrow

  %%%%%%%%%%% hand:
 \put(  12., -12){$\b{t}_{01}=\b{x}_1-\b{x}_0$ (tangent vector) }
\put( -15., -22.){edge-bubble:
     \ $b=\lambda_0\lambda_1p\b{t}_{01}^T\b{t}_{01}$, \
    $b\cdot{\b{n}_i}=0$, \ $i=0,1,2,3$.  }
\put(-15, 38){ Tetrahedron $K$: }
 \multiput(  65,   0)(  -0.9,   0.3){8}{\la}

 \end{picture}\end{center}

\caption{\lab{e-b} An edge-bubble function $b=\lambda_0\lambda_1p\b{t}_{01}^T\b{t}_{01}$,
  $p\in P_{k-2}(K)$, on an edge $\b x_0\b x_1$ of tetrahedron $K$. }
  \end{figure}
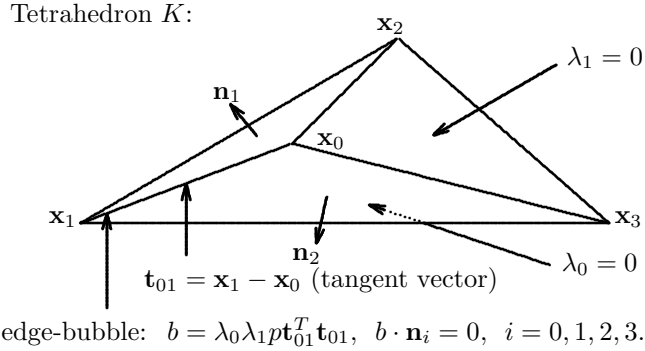

Let $\b x_0$, $\b x_1$, $\b x_2$ and $\b x_3$ be the four vertices of a
   tetrahedron $K$, cf. Figure \mref{e-b}.
The referencing mapping is then
\a{  \b x &= F_K(\hat {\b x})
     = \b x_0 + \p{ \b x_1 -\b x_0 &
                     \b x_2 -\b x_0&
                     \b x_3 -\b x_0  } \hat{\b x}, }
      mapping the reference tetrahedron $\hat K=\{ 0 \le \hat x_1, \hat x_2, \hat x_3,
          1-\hat x_1 -\hat x_2 -\hat x_3 \le 1 \}$ to $K$.
   Then the inverse mapping is
\an{ \lab{imap}  \hat {\b x} &= \p{\b n_1^T\\ \b n_2^T \\ \b n_3^T} (\b x-\b x_0), }
where
 \an{\lab{nnt} \p{\b n_1^T\\ \b n_2^T \\ \b n_3^T }=\p{ \b x_1 -\b x_0 &
                     \b x_2 -\b x_0  &
                     \b x_3 -\b x_0}^{-1}. }
By \meqref{imap}, these normal vectors
        are coefficients of the barycentric variables:
  \a{ \lambda_1 &=  \b n_1 \cdot (\b x-\b x_0), \\
     \lambda_2 &= \b n_2 \cdot (\b x-\b x_0), \\
      \lambda_3 &= \b n_3 \cdot (\b x-\b x_0), \\
    \lambda_0 &= 1-\lambda_1-\lambda_2-\lambda_3. }
On each face triangle, say $\b x_0\b x_2\b x_3$,  all three edges (the tangent vector),
     $ {\b x_0\b x_2}$,  $ {\b x_0\b x_3}$ and $ {\b x_2\b x_3}$,
     are orthogonal to the face normal vector $\b n_1$.
For convenience,  we introduce the tangent vectors and
    their tensors:
   \an{\lab{T} \ad{
     \b t_{01}&=\b x_1 -\b x_0, & T_{01} & = \b t_{01}^T \b t_{01}, \\
       \b t_{02}&=\b x_2 -\b x_0, & T_{02} & = \b t_{02}^T \b t_{02}, \\
       \b t_{03}&=\b x_3 -\b x_0, & T_{03} & = \b t_{03}^T \b t_{03}, \\
       \b t_{12}&=\b x_2 -\b x_1, & T_{12} & = \b t_{12}^T \b t_{12}, \\
       \b t_{23}&=\b x_3 -\b x_2, & T_{23} & = \b t_{23}^T \b t_{23}, \\
       \b t_{13}&=\b x_3 -\b x_1, & T_{13} & = \b t_{13}^T \b t_{13}. } }
With them,  we define the $H(\d,K, {\mathbb S})$ bubble functions
\an{\lab{bh}  \Sigma_{K,b} = \operatorname{span}
    \{ & \lambda_0\lambda_1p_1T_{01},  \lambda_0\lambda_2p_2T_{02},
     \lambda_0\lambda_3p_3T_{03}, \\
      \nonumber & \lambda_1\lambda_2p_4T_{12}, \lambda_2\lambda_3p_5T_{23},
        \lambda_1\lambda_3p_6T_{13}   \},              }
 where $p_{1}, \dots, p_{6}$   are 3D $P_{k-2}$ polynomials.
 Note that each bubble function, say, $\lambda_0\lambda_1p_1T_{01}$,
   vanishes on two face triangles ($\lambda_0=0$, $\lambda_1=0$)
   and has zero normal component on the other two face triangles
    ($T_{01}\cdot \b n_2=\b 0$, $T_{01}\cdot \b n_3=\b 0$.)
 Thus, the matching of $\d \tau_h $ and $ v_h$ is done locally on $K$, independently
   of the matching on neighboring elements. To characterize the bubble space $\Sigma_{K,b}$,  we need the following
   lemma.

   \begin{lemma}\lab{6T}
   The six symmetric tensors $T_{ij}$ in \meqref{T} are linearly independent, and
    form a basis of $\mathbb{S}$.
\end{lemma}

\begin{proof}  Each tensor $T_{ij}=\b t_{ij}  \b t_{ij}^T$ is a positive
  semi-definite matrix, on a tetrahedron $K$.
   We would show that the constants $c_{ij}$ are all
      equal to zero in
  \a{  T= c_{01} T_{01} + c_{02} T_{02} +c_{03} T_{03}
       + c_{12} T_{12} + c_{23} T_{23} +c_{13} T_{13} =0. }
  First, we compute the bilinear form, cf. Figure \mref{e-b}, by \meqref{nnt},
   \a{ \b n_1^T T \b n_1 = c_{01} 1\cdot 1 + c_{02} 0 +c_{03} 0
       + c_{12} (-1)(-1)  + c_{23} 0 +c_{13}(-1)(-1)=0. }
   Here, by \meqref{nnt} and \meqref{T},
   \a{ \b t_{01}^T \b n_1 &=1, \\
       \b t_{12}^T \b n_1 &= (\b t_{02}^T - \b t_{01}^T)\b n_1 =0- 1, \\
       \b t_{13}^T \b n_1 &= (\b t_{03}^T - \b t_{01}^T)\b n_1 =0- 1. }
    Symmetrically, by evaluating $\b n_i T\b n_i$ for $i=0,1,2,3$, where
   $\b n_0=-\b n_1-\b n_2-\b n_3$,
   we have
   \an{ \lab{c-1} \left\{ \ad{
        c_{01}+ c_{02}  + c_{03}  &=0, \\
        c_{01}+ c_{12}  + c_{13}  &=0, \\
        c_{02}+ c_{12}  + c_{23}  &=0, \\
        c_{03}+ c_{13}  + c_{23}  &=0. }
      \right. }
    Note that $\b n_0\ne \b 0$ as $K$ is a non-singular tetrahedron.
    Next, we introduce three (non-unit) vectors $\b s_i$
      orthogonal to the three pairs of skew edges,
      $\overline{\b x_0\b x_1}$ and $\overline{\b x_2\b x_3}$,
      $\overline{\b x_0\b x_2}$ and $\overline{\b x_1\b x_3}$,
      $\overline{\b x_0\b x_3}$ and $\overline{\b x_1\b x_2}$, respectively,
        cf. Figure \mref{e-b}.
      That is,  \a{ \b s_1=  \frac{ \b t_{01} \times \b t_{23} }
                {6|K|}, }
            because $|K|\ne 0$ and consequently $|\b t_{01} \times \b t_{23} |\ne0$.
       Thus $\b s_1\cdot \b t_{01}=0$, $\b s_1\cdot \b t_{02}=-1$,
               $\b s_1\cdot \b t_{03}=-1$, $\b s_1\cdot \b t_{12}=-1$,
               $\b s_1\cdot \b t_{13}=-1$, and $\b s_1\cdot \b t_{23}=0$.
      By evaluating $\b s_i^T T \b s_i$,  it follows that
   \an{ \lab{c-2} \left\{ \ad{
        c_{02}+ c_{03}  + c_{12}   + c_{13}  &=0, \\
        c_{01}+ c_{03}  +  c_{12}  + c_{23}  &=0, \\
        c_{01}+ c_{02}  + c_{13}  + c_{23}  &=0. }
      \right. }
       By the first two equations in \meqref{c-1} and the first equation in
         \meqref{c-2}, we get
      \a{ 2 c_{01} =0 .}
      Symmetrically,  we find all $c_{ij}=0$.
      Thus $\{ T_{ij}\}$ is a linearly independent set of tensors.
      As $\dim \mathbb{S}=6$,   $\{ T_{ij}\}$ is a basis.
\end{proof}

An equivalent but more practical definition of the stress finite element
   space $\Sigma_h$ is
\an{ \lab{Sh2} \Sigma_h= \Big\{ \sigma=\sigma_a+\sigma_b
     \in H(\d,\Omega,\mathbb{S}) \ \Big| \
     \sigma_a \in \widetilde \Sigma_h,\ \sigma_b|_K\in \Sigma_{K,b}
     \   \forall K\in\mathcal{T}_h \Big\}, }
where $\widetilde\Sigma_h$ and $\Sigma_{K,b}$ are defined in
    \meqref{tSh} and \meqref{bh},  respectively.

 It follows from the definition of $V_h$ ($P_{k-1}$ polynomials)
    and $\Sigma_h$ ($P_k$ polynomials) that
   \a{ \d  \Sigma_h \subset V_h.}
This, in turn, leads to a strong divergence-free space:
 \an{ \lab {kernel}
    Z_h&:= \{\tau_h\in\Sigma_h \ | \ (\d\tau_h, v)=0 \quad
	\hbox{for all } v\in V_h\}\\
    \nonumber
          &= \{\tau_h \in\Sigma_h \ | \  \d \tau_h=0
    \hbox{\ pointwise } \}.
    }

The mixed finite element approximation of Problem (1.1) reads: Find
   $(\sigma_h,~u_h)\in\Sigma_h\times V_h$ such that
 \e{ \left\{ \ad{
    (A\sigma_h, \tau)+({\rm div}\tau, u_h)&= 0 &&
              \hbox{for all \ } \tau \in\Sigma_h,\\
     (\d\sigma_h, v)& = (f, v) &&  \hbox{for all \ } v\in V_h.
      } \right. \lab{DP}
    }

\section{Stability and convergence}
The convergence of the finite element solutions follows
   the stability and the standard approximation property.
So we consider first the well-posedness  of the discrete problem
    \meqref{DP}.
By the standard theory,  we only need to prove
   the following two conditions, based on their counterpart at
    the continuous level.

\begin{enumerate}
\item K-ellipticity. There exists a constant $C>0$, independent of the
   meshsize $h$ such that
    \an{ \lab{below} (A\tau, \tau)\geq C\|\tau\|_{H(\d)}^2\quad
       \hbox{for all } \tau \in Z_h, }
    where $Z_h$ is the divergence-free space defined in \meqref{kernel}.

\item  Discrete B-B condition.
    There exists a positive constant $C>0$
            independent of the meshsize $h$, such that
    \an{\lab{inf-sup}
   \inf_{0\neq v\in V_h}   \sup_{0\neq\tau\in\Sigma_h}\frac{({\rm
        div}\tau, v)}{\|\tau\|_{H(\d)}  \|v\|_{L^2(\Omega)} }\geq
    C .}
\end{enumerate}

It follows from $\d  \Sigma_h \subset V_h$ that $\d  \tau=0$ for
   any $\tau\in Z_h$. This implies the above K-ellipticity condition
	\meqref{below}.
It remains to show the discrete B-B condition \meqref{inf-sup},
  in the following two lemmas.

\begin{lemma}\label{lemma1}
For any $v_h\in V_h$,  there is a $\tau_h \in
    \widetilde \Sigma_h\subset  \Sigma_h$  such that,
  for all polynomial $p\in P_{k-3}(K,\mathbb{R}^3)$, $K\in\mathcal{T}_h$,
   \bq\label{l-1}  \int_K (\d\tau_h-v_h) \cdot p\, d\b x=0
      \quad \hbox{\rm
      and } \quad \|\tau_h\|_{H(\d)}\leq C\|v_h\|_{L^2(\Omega)}. \eq
     \end{lemma}

\begin{proof}
  Let $v_h\in V_h$.
  By the stability of the continuous formulation,
      cf. \cite{Arnold-Winther-conforming}, there is
    a $\tau \in  H^1(\Omega,\mathbb{S})$ such that,
   \a{ \d\tau=v_h \quad \hbox{\rm
      and } \quad \|\tau\|_{H^1(\Omega)}\le  C\|v_h\|_{L^2(\Omega)}. }
As $\tau \in H^1(\Omega,\mathbb{S})$,  we modify
   the Scott-Zhang \cite{Scott-Zhang}
   interpolation operator slightly to define a flux preserving
   interpolation.
 \a{  I_h \ : \ H^1(\Omega,\mathbb{S})
              & \to \Sigma_h \cap H^1(\Omega,\mathbb{S}) =\widetilde \Sigma_h , \\
         \ \tau & \mapsto \tau_h =: I_h \tau. }
Here the interpolation is done inside a subspace,
   the continuous finite element subspace
     $\Sigma_h \cap H^1(\Omega,\mathbb{S}) $.
$I_h \tau$ is defined by its values at the Lagrange nodes.

At a vertex node or a node inside an edge, $\b x_i$,
   $I_h\tau(\b x_i)$ is defined as the nodal value of $\tau$ at the
     point if $\tau$ is continuous, but in general, $I_h\tau(\b x_i)$ is
       defined as an average value on a face triangle,  on whose edge the node
          is,
        as in \cite{Scott-Zhang}.
After defining the nodal values at edges of tetrahedra,
  the nodal values of $\tau_h$
      at the nodes inside each face triangle $T$ of a tetrahedron
          are defined by
   the $L^2$-orthogonal projection on the triangle $T$:
  \an{ \lab{a-1}
    \int_T \tau_{h,ij}  p  \, dS &=
    \int_T \tau_{ ij}  p  \, dS\quad
   \forall p\in P_{k-3}(T, \mathbb{R}),   }
$i,j=1,2,3$,
 where $\tau_{h,ij}$ and $\tau_{ij}$ are the $(i,j)$-th components of
    $\tau_h$ and $\tau$, respectively,
   and $T$ is a face triangle of a tetrahedron in the tetrahedral
      triangulation $\mathcal{T}_h$.
The number of equations in \meqref{a-1} is same as the number of
   internal degrees of freedom of $P_k$ polynomials,  $\dim P_{k-3}$.
At the Lagrange nodes inside a tetrahedron,
   $I_h\tau(\b x_i)$ is defined by
   the $L^2$-orthogonal projection on the tetrahedron, satisfying
 \an{ \lab{a-2}
    \int_K  \tau_{h,ij}  p  \, d\b x &= \int_K \tau_{ij} p\, d\b x  \quad
   \forall p\in P_{k-4}(K,\mathbb{R}),   }
    where $K$ is an element of  $\mathcal{T}_h$.
  It follows by the stability of the Scott-Zhang operator that
   \a{ \| I_h \tau \|_{H^1(\Omega)}\le C\| \tau \|_{H^1(\Omega)}\le C\|v_h\|_{L^2(\Omega)}. }
In particular,
  \a{ \|I_h \tau\|_{H(\d)} \le \| I_h \tau \|_{H^1(\Omega)} \le C\|v_h\|_{L^2(\Omega)}. }
By \meqref{a-1} and \meqref{a-2},
  we get  a partial-divergence matching property of $I_h$:
     for any $p\in P_{k-3}(K,\mathbb{R}^3)$,
    as the symmetric gradient $\epsilon( p)\in P_{k-4}(K,\mathbb{S})$,
 \a{  \int_K (\d\tau_h-v_h) \cdot p \, d\b x
     & = \int_{\partial K}  (\tau_h  \b n ) \cdot p\,ds
        -\int_K \tau_h : \epsilon( p) \, d \b x
           - \int_K v_h \cdot p\, d\b x \\
        & = \int_{\partial K}  (\tau  \b n)\cdot p\,  ds
        -\int_K \tau: \epsilon( p ) \, d \b x
           - \int_K v_h \cdot p\,d\b x \\
	 	 &= \int_K (\d\tau -v_h) \cdot p \, d\b x = 0. }
\end{proof}

\begin{lemma}\lab{lemma2}
For any $v_h\in V_h$, if
  \an{ \lab{o3} \int_K  v_h \cdot p \, d\b x=0
     \quad\hbox{ for all $p\in P_{k-3} (K,\mathbb{R}^3)$ and
   all $K\in\mathcal{T}_h$, } }
   there is a $\tau_h \in \Sigma_h$  such that
   \an{ \lab{l-2}  \d\tau_h = v_h
      \quad \hbox{\rm
      and } \quad \|\tau_h \|_{H(\d)}\leq C\|v_h\|_{L^2(\Omega)}. }
    \end{lemma}

\begin{proof}
As we assume polynomial degree $k\ge 4$ in $V_h$,
  $p\in P_{k-3} (K,\mathbb{R}^3)\supset  P_1 (K,\mathbb{R}^3)\supset
     R(K)$ where $R(K)$ is the set of 6-dimensional, local rigid motions:
  \an{\lab{R}
    R(K) = \Big\{
     \p{ a_1 - a_4 y -a_5 z \\
         a_2 + a_4 x -a_6 z \\
         a_3 + a_5 x +a_6 y } \ | \ a_1, a_2,a_3,a_4,a_5,a_6\in\mathbb{R}
        \Big\}. }
 So if $v_h$ satisfies \meqref{o3}, $v_h$ is in the following
    local rigid-motion free space:
\an{ \lab{Vo} V_{h,\perp R}
      =\Big\{ v_h\in V_h \ | \
            \int_K v_h \cdot p \,d\b x = 0 \ \forall p\in R(K) \ \hbox{and }
              \forall K\in \mathcal{T}_h\Big\}. }
We will prove a stronger result that if $v_h\in V_{h,\perp R}$,
    then there is a $\tau_h$ satisfying \meqref{l-2}.
This $\tau_h$ is constructed, according to $v_h$, on each element $K$,
   independently of the construction on neighboring elements.
On one element $K$,  we show $\d \Sigma_{K,b}=V_{h,\perp R}|_K$ where
    $\Sigma_{K,b}$ is the edge-bubble  space, defined
    in \meqref{bh}.
If $\d \Sigma_{K,b} \ne V_{h,\perp R}|_K$,
    there is a nonzero $ v_h\in V_{h,\perp R}$
    such that
  \a{ \int_K \d \tau_h \cdot v_h\, d\b x=0 \quad \forall \tau_h \in \Sigma_{K,b}. }
By integration by parts,   for  $\tau_h \in \Sigma_{K,b}$,
  we have
  \an{ \lab{o5}
   \int_K  \d \tau_h \cdot
     v_h d\b x = \int_{ K} \tau_h: \epsilon(v_h) d\b x=0, }
  where $\epsilon(v_h)$ is the symmetric gradient, $(\nabla v_h + \nabla^T v_h)/2$.

Let $\{M_{ij}, i=0,1,2, \, j=i,\dots,3,\}$ be the dual basis of
   the symmetric space,  under $\mathbb R^9$ inner-product, of $\{T_{ij}\}$,
   defined in \meqref{T},
i.e.
  \an{\lab{dual} M_{ij}=M_{ij}^T,  \quad  M_{ij} \cdot T_{i'j'}
    =\delta_{ij,i'j'}. }
For example, if $K$ is the unit right tetrahedron,  then $\{T_{ij}\}$ would be
 \a{ & & &\p{1&0&0\\ 0&0&0\\  0&0&0}, &&\p{0&0&0\\ 0&1&0\\  0&0&0},
         &&\p{0&0&0\\ 0&0&0\\  0&0&1},& & \\
      & & &\p{1&-1&0\\ -1&1&0\\  0&0&0}, &&\p{1&0&-1\\ 0&0&0\\ -1&0&1},
         &&\p{0&0&0\\ 0&1&-1\\  0&-1&1},& & }
   and the unique $\{M_{ij}\}$ would be
 \a{ & & &\p{1&1/2&1/2\\1/2&0&0\\  1/2&0&0},
        &&\p{0&-1/2&0\\-1/2&1&1/2\\  0&1/2&0},
         &&\p{0&0&1/2\\ 0&0&1/2\\  1/2&1/2&1},& & \\
      & & &\p{0&-1/2&0\\  -1/2&0&0\\  0&0&0},
      &&\p{0&0&-1/2\\ 0&0&0\\ -1/2&0&0},
         &&\p{0&0&0\\ 0&0&-1/2\\  0&-1/2&0}.& & }

Under the dual basis, we have a unique expansion, as $\epsilon(v_h)
    \in  P_{k-2}(K,\mathbb{S})$,
  \an{\lab{expan} \epsilon(v_h)= q_1 M_{01} + q_2 M_{02}+ q_3M_{03}
              q_4 M_{12} + q_5 M_{23}+ q_6M_{13} , }
               for some  $q_i \in P_{k-2}(K)$.
Selecting $\tau_1 = \lambda_0\lambda_1q_1 T_{01}\in \Sigma_{K,b}$,
         we have, by \meqref{dual},
   \a{ 0 = \int_K \tau_1: \epsilon(v_h) d\b x
         = \int_K \lambda_0\lambda_1 q_1^2(\b x) d\b x. }
As $\lambda_0\lambda_1>0$ on $K$,  we conclude that $q_1\equiv 0$.
Similarly, the other five $q_i$ in \meqref{expan} are zero.
Thus, by \meqref{o5}, $v_h\equiv 0$ and $\d \sigma_{K,b}=V_{h,\perp R}|_K$.
As the matching $\d \tau_h=v_h$ is done on one element $K$,
   by affine mapping and scaling argument,  \meqref{l-2} holds.

\end{proof}

 We are in the position to show the well-posedness of the discrete problem.
\begin{lemma}
 For the discrete problem (\ref{DP}), the K-ellipticity \meqref{below}
    and the discrete B-B
 condition \meqref{inf-sup} hold uniformly.
  Consequently,  the discrete
     mixed problem \meqref{DP} has a unique solution
         $(\sigma_h,~u_h)\in\Sigma_h\times V_h$.
\end{lemma}
\begin{proof}  The  K-ellipticity immediately follows from the fact
      that $\d  \Sigma_h \subset V_h$.
      To prove the  discrete B-B  condition \meqref{inf-sup},
      for any $v_h\in V_h$,
        it follows from Lemma \ref{lemma1} that there exists a
      $\tau_{1}\in \Sigma_h$ such that,  for any polynomial $p\in P_{k-3}(K, \mathbb{R}^3)$,
   \bq  \int_K (\d\tau_1-v_h) \cdot pd\b x=0
      \quad \hbox{\rm
      and } \quad \|\tau_1\|_{H(\d)}\leq C\|v_h\|_{L^2(\Omega)}. \eq
Then it follows from Lemma \ref{lemma2} that
 there is a $\tau_2 \in \Sigma_h$  such that
   \an{   \d\tau_2 = v_h-\d\tau_1
      \quad \hbox{\rm
      and } \quad \|\tau_2 \|_{H(\d)}\leq C\|\d\tau_1-v_h\|_{L^2(\Omega)},}
Let $\tau=\tau_1+\tau_2$.  This implies that
\begin{equation}
\d\tau=v_h \text{ and } \|\tau\|_{H(\d)}\leq C\|v_h\|_{L^2(\Omega)},
\end{equation}
this proves the discrete B-B condition \meqref{inf-sup}.
\end{proof}

\begin{theorem}\label{MainError} Let
  $(\sigma, u)\in\Sigma\times V$ be the exact solution of
   problem \meqref{eqn1} and $(\tau_h, u_h)\in\Sigma_h\times
   V_h$ the finite element solution of \meqref{DP}.  Then,
   for $k\ge 4$,
\an{ \lab{t1} \|\sigma-\sigma_h\|_{H({\rm div})}
    + \|u-u_h\|_{L^2(\Omega)}&\le     Ch^k(\|\sigma\|_{H^{k+1}(\Omega)}+\|u\|_{H^k(\Omega)}).
      }
\end{theorem}

\begin{proof}
 The stability of the elements and the standard theory of mixed
  finite element methods \cite{Brezzi, Brezzi-Fortin} give the
  following quasioptimal error estimate immediately
\an{
  \label{theorem-err1} \|\sigma-\sigma_h\|_{H({\rm
  div})}+\|u-u_h\|_{L^2(\Omega)}\leq C \inf\limits_{\tau_h\in\Sigma_h,v_h\in
  V_h}\left(\|\sigma-\tau_h\|_{H({\rm div})}+\|u-v_h\|_{L^2(\Omega)}\right).}
Let $P_h$ denote the local $L^2$ projection operator,
   or triangle-wise interpolation operator,  from $V$ to $V_h$,
  satisfying the error estimate
\an{\label{proj-error}
   \|v-P_hv\|_{L^2(\Omega)}\leq Ch^k\|v\|_{H^k(\Omega)} \text{ for any }v\in H^k(\Omega, \mathbb{R}^3). }
Choosing $\tau_h=I_h\sigma\in \Sigma_h$
    where $I_h$ is defined in \meqref{a-1} and \meqref{a-2},
we have \cite{Scott-Zhang},  as $I_h$ preserves symmetric $P_k$ functions locally,
   \an{ \lab{p-err2}
       \|\sigma -\tau_h\|_{L^2(\Omega)} + h |\sigma -\tau_h|_{H(\d)}
        \le Ch^{k+1} \|\sigma\|_{H^{k+1}(\Omega)}. }
 Let $v_h= P_h v$ and $\tau_h=I_h\sigma$ in (\ref{theorem-err1}),
  by (\ref{proj-error}) and \meqref{p-err2}, we
   obtain  \meqref{t1}.
\end{proof}

\section{Numerical tests}\lab{s-numerical}

We compute one example in 3D, by $P_4$ and by $P_5$ mixed finite element methods.
It is a pure displacement problem  on the unit cube
   $\Omega=(0,1)^3$ with a homogeneous boundary condition
         that $u\equiv 0$ on $\partial\Omega$.
In the computation, we let
   \a{
      A \sigma &= \frac 1{2\mu} \left(
       \sigma - \frac{\lambda}{2\mu + n \lambda} \operatorname{tr}(\sigma)
        \delta \right), \quad n=3, }
  where $\delta=\p{1 &0&0\\0&1&0\\0&0&1}$, and $\mu=1/2$ and $\lambda=1$ are the
    Lam\'e constants.

Let  the  exact solution on the unit square $[0,1]^3$ be
   \e{\lab{e1}   u= \p{2^4 \\2^5 \\ 2^6 } x(1-x)y(1-y)z(1-z). }
Then, the true stress function $\sigma$
     and the load function $f$ are defined by the equations in
    \meqref{eqn1},   for the given  solution $u$.

\begin{figure}[htb] \setlength\unitlength{0.8pt}
  \begin{center}  \begin{picture}(320,100)(0,0)
 \put(0,0){\begin{picture}(100,100)(0,0)
     \multiput(0,0)(80,0){2}{\line(0,1){80}}     \multiput(0,0)(0,80){2}{\line(1,0){80}}
    \multiput(80,0)(0,80){2}{\line(2,1){40}}    \multiput(0,80)(80,0){2}{\line(2,1){40}}
    \multiput(80,0)(40,20){2}{\line(0,1){80}}    \multiput(0,80)(40,20){2}{\line(1,0){80}}
    \multiput(0,0)(40,0){1}{\line(1,1){80}}
    \multiput(80,0)(0,80){1}{\line(2,5){40}}    \multiput(0,80)(80,0){1}{\line(6,1){120}}
   \end{picture} }
 \put(180,0){\begin{picture}(100,100)(0,0)
     \multiput(0,0)(40,0){3}{\line(0,1){80}}     \multiput(0,0)(0,40){3}{\line(1,0){80}}
    \multiput(80,0)(0,40){3}{\line(2,1){20}}    \multiput(0,80)(40,0){3}{\line(2,1){20}}
    \multiput(100,10)(0,40){3}{\line(2,1){20}}    \multiput(20,90)(40,0){3}{\line(2,1){20}}
    \multiput(80,0)(20,10){3}{\line(0,1){80}}    \multiput(0,80)(20,10){3}{\line(1,0){80}}
    \multiput(0,0)(40,0){1}{\line(1,1){80}}     \multiput(40,0)(-40,40){2}{\line(1,1){40}}

    \multiput(80,0)(0,80){1}{\line(2,5){40}}    \multiput(0,80)(80,0){1}{\line(6,1){120}}
    \multiput(100,10)(-20,30){2}{\line(2,5){20}} \multiput(40,80)(-20,10){2}{\line(6,1){60}}
   \end{picture} }
    \end{picture} \end{center}
\caption{ \lab{grid} The initial grid for
 \meqref{e1}, and its level 2 refinement. }
\end{figure}
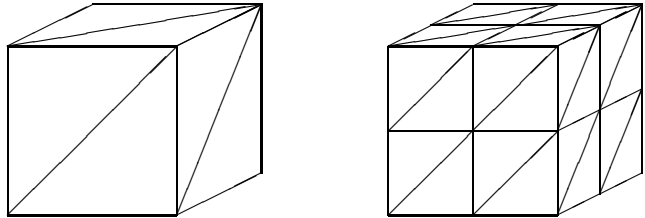

In the computation, the level one grid is the given domain with a diagonal
   line shown in Figure \ref{grid}.
Each grid is refined into a half-sized grid uniformly, to get
   a higher level grid, shown in Figure \ref{grid}.
In all the computation, the discrete systems of equations are
  solved by Matlab backslash solver.
In Table \mref{b1}, the errors and the convergence order
   in various norms are listed  for the true solution \meqref{e1},
   by $P_4$ mixed finite elements in \meqref{Sh} and \meqref{Vh}, with
    $k=4$ there.
Here  $I_h$ is the usual nodal
   interpolation operator.
The optimal order of covergence is achived in Table \mref{b1},
   confirming Theorem \mref{MainError}.

\begin{table}[htb]
  \caption{ The error and the order of convergence by $P_4$ finite elements,
    $k=4$ in \eqref{Vh} and  \meqref{Sh}, for \meqref{e1}.}
\lab{b1}
\begin{center}  \begin{tabular}{c|cc|cc|cc}  %\multispan{3}
\hline &  $ \|I_h\sigma -\sigma_h\|_{L^2(\Omega)}$ & $h^n$  & $ \|I_h u- u_h\|_{L^2(\Omega)}$ &$h^n$ &
    $ \|\d(I_h\sigma -\sigma_h)\|_{L^2(\Omega)}$ & $h^n$  \\ \hline
 1&    0.33567012&0.0&    0.05860521&0.0&    3.41111411&0.0\\
 2&    0.02041247&4.0&    0.00661542&3.1&    0.21319463&4.0\\
 3&    0.00125425&4.0&    0.00044841&3.9&    0.01332466&4.0\\      \hline
\end{tabular}\end{center} \end{table}

In Table \mref{b2}, the errors and the convergence order
   in various norms are listed  for the true solution \meqref{e1},
   by $P_5$ mixed finite elements in \meqref{Sh} and \meqref{Vh}, with
    $k=5$ there.
Here the exact solution $\sigma$ is a polynomial tensor of degree 5.
Thus,  it is in the stress finite element space $\Sigma_h$ and
   the finite element solution $\sigma_h$ is exact.
It is computed so,  shown in the second column and the six column
   in Table \ref{b2}.
The optimal order of covergence is achived for the displacement $u$
   in Table \mref{b2}  (up to the computer accuracy),
   confirming Theorem \mref{MainError}.

\begin{table}[htb]
  \caption{ The error and the order of convergence by $P_5$ finite elements,
    $k=5$ in \eqref{Vh} and  \meqref{Sh}, for \meqref{e1}.}
\lab{b2}
\begin{center}  \begin{tabular}{c|cc|cc|cc}  %\multispan{3}
\hline &  $ \|I_h\sigma -\sigma_h\|_{L^2(\Omega)}$ & $h^n$  & $ \|I_h u- u_h\|_{L^2(\Omega)}$ &$h^n$ &
    $ \|\d(I_h\sigma -\sigma_h)\|_{L^2(\Omega)}$ & $h^n$   \\ \hline
 1&    0.00000002&0.0&    0.01937914&0.0&    0.00000011&0.0\\
 2&    0.00000002&0.0&    0.00089726&4.4&    0.00000031&0.0\\ \hline
\end{tabular}\end{center} \end{table}

\end{document}